\documentclass[11pt]{amsart}
\usepackage{amsmath,amsfonts,amsthm,amscd,amssymb,mathrsfs,amssymb,bm}
\usepackage{mathrsfs}
\usepackage[dvips]{graphicx}
\usepackage{wrapfig}
\usepackage[all]{xy}
\newtheorem{theorem}{Theorem}

\newtheorem{proposition}[theorem]{Proposition}
\newtheorem{lemma}[theorem]{Lemma}
\newtheorem{definition}[theorem]{Definition}

\newtheorem{remark}[theorem]{Remark}
\newcommand{\CP}{\mathbb{CP}}
\newcommand{\CC}{\mathbb{C}}

\newcommand{\RR}{\mathbb{R}}
\newcommand{\ZZ}{\mathbb{Z}}

\newcommand{\ol}{\overline}
\newcommand{\lra}{\longrightarrow}

\newcommand{\set}{\,|\,}
\newcommand{\proofend}{\hfill$\square$}

\newcommand{\Bs}{{\rm{Bs}}}

\newcommand{\vsp}{\vspace{3mm}}

\newcommand{\aaa}{{\alpha}}

\numberwithin{equation}{section}
\numberwithin{theorem}{section}
\begin{document}
\bibliographystyle{alpha} 
\title[Geometry of generic Moishezon twistor spaces
on $4\mathbb{CP}^2$ II]{Geometry of generic Moishezon twistor spaces
on $4\mathbb{CP}^2$ II:  degenerate cases}
\author{Nobuhiro Honda}
\address{Mathematical Institute, Tohoku University,
Sendai, Miyagi, Japan}
\email{honda@math.tohoku.ac.jp}
\begin{abstract}
We continue to study  twistor spaces on 
the connected sum of four complex projective planes, whose
anticanonical map is of degree two over the image.
In particular, we determine the defining equation of the
branch divisor of the anticanonical map in an explicit form.
Together with previous two articles (\cite{HonDSn1} and \cite{HonDS4_1}),
this completes explicit description of all such 
twistor spaces.
\end{abstract}
\maketitle
\setcounter{tocdepth}{1}
\vspace{-5mm}

\section{Introduction}
This paper is a sequel to an article \cite{HonDS4_1},
where we intensively studied twistor
spaces on $4\CP^2$ which could be considered as the
most generic ones among all Moishezon twistor spaces on $4\CP^2$.
The most characteristic property of these twistor spaces
is that they have a natural structure of double covering
over a very simple rational threefold, so called a scroll
of planes over a conic.
Moreover, the branch divisor of the covering is always a 
cut of the scroll by a {\em quartic hypersurface.}
This allowed us to describe each of these twistor spaces by
a single quartic polynomial.
It was shown that this quartic polynomial takes
not general but a very special form.
This description has a clear advantage
when describing a global structure of the moduli space of these twistor spaces.

While these twistor spaces are most generic  among
all Moishezon twistor spaces on $4\CP^2$ as above,
we showed in another article \cite{HonDS4_3}
(in which we classified all Moishezon twistor spaces
on $4\CP^2$)
that there still exist  other three families of twistor spaces
on $4\CP^2$ which have a natural structure of 
a double covering over the same scroll.
The twistor spaces in these families are strictly different from the ones in \cite{HonDS4_1},
but they are obtained from the generic ones 
by taking a limit under deformations.
In particular, the branch divisor of the double covering is
still a cut of the scroll by a quartic hypersurface, while its defining equation should be
subject to  stronger  constraint.

Among the three families, the most special one is nothing but the family of twistor spaces studied in \cite{HonDSn1} (specialized to the case of $4\CP^2$).
In particular for twistor spaces in this family the defining equation of the quartic hypersurface
is already determined.
The main purpose of this article is to obtain the equation of the quartic hypersurface for the remaining two families
of the spaces

As explained above, twistor spaces whose 
anticanonical map is of degree two can be classified 
into four types.
Distinction of the types can be detected
by the number of irreducible components of the base curve
for the half-anticanonical system on the twistor spaces.
We just call these as type I, II, III and IV
(Definition \ref{def:type};
this naming is justified when we obtain defining equation
of the quartic hypersurface.).
Type I is most generic (so they are studied in \cite{HonDS4_1})
and type IV is most special (so they are studied in 
\cite{HonDSn1}).
Hence in this paper we study those of type II and type III.
In Section \ref{s:acs} we take a general member of the 
half-anticanonical system and 
investigate  structure of its bi-anticanonical
system.
The last system exhibits the surface as a double covering
of $\CP^2$ with branch being a quartic curve,
whose singularities depend on the types I--IV
(as displayed in Figure \ref{fig:quartics}). 

As above, the structure of the twistor spaces
can be captured through the anticanonical system.
In Section 3, we find {\em reducible} anti-canonical divisor(s) of the twistor space, which will be a key
for determining the equation of the branch divisor 
of the anticanonical map.

In contrast with the generic ones studied in \cite{HonDS4_1}, the base locus of the anticanonical map
of the present twistor spaces is quite complicated.
In Section 4 we give a full elimination of this locus through explicit blowups.
In particular, we find that the anticanonical map
contracts some divisors to curves.
(For the case of type I such a divisor does not exist
as investigated in \cite{HonDS4_1}.)
This information will give a strong constraint for the 
equation of the branch divisor of the anticanonical map.

In Section 5, by using the reducible anticanonical 
divisor(s) found in 
Section 3, we first find five special hyperplanes in $\CP^4$
whose intersection with the scroll touches the branch divisor along a curve.
These curves are called double curves.
Next we show that there  exists a hyperquadric
in $\CP^4$ which contains all these five double curves.
Finally by using the hyperquadric and also the information about the anticanonical map obtained in 
Section \ref{s:analysis}, we determine
defining equation of the branch divisor
of the anticanonical map, for both cases of type II and type III
(Theorem \ref{thm:main}).
The argument in the proof is mostly algebraic.
We also give an account about how
the defining equation of the quartic hypersurface
degenerates when the twistor space changes the type.

In Section 6 we compute dimension of the moduli spaces of these twistor spaces.
The conclusion including the cases of type I and type IV
is as follows:

\vsp
\begin{center}
\begin{tabular}{|c||c|c|c|c|}
\hline
                & type I & type II & type III & type IV\\ \hline
dimension of the moduli space  & 9 &   7  &  5    & 4  \\ 
dimension of the automorphism group & 0      &    0    &   0      & 1\\ \hline
\end{tabular} 
\end{center}

\vsp\noindent
Thus whole  picture is now understood to a
considerable degree.

Finally, as mentioned in \cite{HonDS4_1}, it looks quite certain that 
the twistor spaces which have the structure of the double covering over
the scroll can be generalized to the ones over $n\CP^2$, $n$ being arbitrary.
The results in \cite{HonDSn1} mean this is  actually the case for type IV spaces.
However, after writing \cite{HonDS4_1}, 
the author noticed that 
{\em the twistor spaces on $4\CP^2$ studied in \cite{HonDS4_1} 
(i.e.\,type I spaces) cannot be generalized to $n\CP^2$,
as long as we stick to the linear system $|(n-2)K^{-1/2}|$.}
On the other hand, for those of type II and type III, 
there seems to be a good chance for generalization by using 
$|(n-2)K^{-1/2}|$, like \cite{HonDSn1}.
This is a reason why we study these cases in depth.

\vsp
\noindent
{\bf Notations.}
For a twistor space, the natural square root of the 
anticanonical bundle is denoted by $F$.
(Hence $2F$ is the anticanonical bundle.) 
The dimension of a linear system always means
the dimension of the parameter projective space.
For a line bundle $L\to X$, 
we write $h^i(X,L) = \dim H^i(X,L)$.
For $s\in H^0(X,L)$ with $s\neq 0$, we denote $(s)$ for the zero-divisor
of $s$.
For  a curve $C$ and a divisor $D$ on $X$,
the intersection number of $C$ and $D$ is denoted
by $(C,D)_X$ or just $(C,D)$.

\section{The anticanonical system of the twistor spaces}
\label{s:acs}
We first make it clear which twistor spaces on $4\CP^2$ we are going to investigate.
For this we recall the following result which is one of  a consequence 
from the classification of all Moishezon twistor spaces,
obtained in \cite{HonDS4_3}:

\begin{proposition}\label{prop:classify}
{\em (\cite[Theorem 1.1]{HonDS4_3})}
Let $Z$ be a Moishezon twistor space on $4\CP^2$ and suppose that the anticanonical map $\Phi=\Phi_{|2F|}$ is (rationally) of degree two over the image.
Then we have the following.
(i) $\dim |2F| = 4$ and the image $\Phi(Z)\subset\CP^4$ is  a scroll of planes over a conic,
(ii) $\dim |F| = 1$ and $\Bs \, |F|$ is a cycle of smooth rational curves,
(iii) the cycle 
consists of $4, 6, 8$ or $10$ irreducible components.
\end{proposition}

Here, by {\em the scroll of planes over a conic},
we mean the inverse image $\pi^{-1}(\Lambda)$,
where $\pi:\CP^4\to\CP^2$ is a linear projection
and $\Lambda$ is a conic in $\CP^2$.
Namely, the scroll is a union of all planes in $\CP^4$
which contain a fixed line.
Clearly such a scroll is unique up to projective 
transformations of $\CP^4$.

An obvious relation between $|F|$ and $|2F|$ 
gives the following commutative diagram
\begin{equation}\label{016}
 \CD
 Z@>\Phi_{|2F|}>> \CP^4\\
 @V\Phi_{|F|} VV @VV\pi V\\
 \CP^1@>>> \CP^2.
 \endCD
 \end{equation}
where the bottom arrow is an embedding of
$\CP^1$ onto the conic $\Lambda\subset\CP^2$.

Throughout this paper we denote 
the cycle $\Bs\,|F|$ appeared in (ii) 
of Proposition \ref{prop:classify}
by the letter $C$
(as in \cite{HonDS4_3}).
Of course, this cycle $C$ is itself real. 
The number of its irreducible components  is significant 
because it is directly connected with the structure of the branch divisor
of the degree two rational map $\Phi$ (over the scroll) in (i)
of the proposition.
So we introduce the following

\begin{definition}\label{def:type}{\em 
Let $Z$ be as in Proposition \ref{prop:classify}.
Then according to the number $4, 6, 8$ or $10$ of the irreducible components 
of the 
cycle $C$, 
 we call $Z$ is of {\em type I,  II,  III} or {\em  IV} respectively.
}
\end{definition}

We note that the twistor spaces studied in \cite{HonDS4_1} are exactly
those of type I, and that the  twistor spaces studied in \cite{HonDSn1} are those of type IV
if we substitute $n=4$ in the paper.
In particular in these papers a defining equation of 
the branch divisor on the scroll is explicitly obtained.
So in this article we are concerned with the cases of type II and type III.
We also mention that among these four kinds of the twistor spaces, type I is most generic and type IV is most special
in the sense that, if I $\le A< B\le$ IV, then type A is obtained
from type B by small deformation;
in other words, type B is obtained as a limit through a
family of type A twistor spaces.
So it might be possible to say that the present 
twistor spaces 
are  intermediately degenerate ones among
all twistor spaces (on $4\CP^2$) whose 
anticanonical map is of degree two.

Let $Z$ be as in Proposition \ref{prop:classify} and 
 $S\in |F|$  any real irreducible  member, which is 
always smooth rational surface with $K_S^2=0$ by \cite{PP94}.
As $\dim|F| =1 $, we have $\dim |K_S^{-1}| = 0$, and 
the unique anticanonical curve is exactly the cycle $C$.
By reality we can write it as 
\begin{align}\label{cycle}
C = \sum _{i=1}^k C_i + \sum _{i=1}^k \ol C_i, 
\end{align}
where $k=2,3,4,5$ according to type I, II, III, IV respectively.
Here we are taking the numbering for the components
in a natural way that $C_i$ and $C_{i+1}$ intersect for $1\le i\le k-1$ and 
$C_k$ and $\ol C_1$ intersect.
Then
by \cite{HonDS4_3} the sequence obtained by arranging the self-intersection numbers
 of the components is respectively given as follows (after a proper cyclic permutation
and an exchange of orientation):
\begin{align}\label{I}
-3,-1,-3,-1 \quad&{\text{for type I}},\\
-3,-2,-1,-3,-2,-1 \quad&{\text{for type II}},\label{II}\\
-3,-2,-2,-1,-3,-2,-2,-1 \quad&{\text{for type III}},\label{III}\\
-3,-2,-2,-2,-1,-3,-2,-2,-2,-1 \quad&{\text{for type IV}}.
\end{align}
These indicate that, for example in the case of type II,  the self-intersection numbers
of the components $C_1,C_2,C_3,\ol C_1,\ol C_2,\ol C_3$ in $S$ are
$-3,-2,-1,-3,-2,-1$ respectively.
As we will see in Section \ref{s:analysis}, 
these numbers are directly related with 
birational geometry of the twistor spaces.



From the relation $2F|_S\simeq 2K_S^{-1}$,
the structure of the anticanonical map of the 
twistor spaces may be investigated via
the bi-anticanonical system of the  surface $S$.
For type I (resp.\,type IV) twistor spaces the last system
is investigated in \cite[Section 2]{HonDS4_1} (resp.\,\cite[Section 2.2]{HonDSn1}; substitute $n=4$).
We now write down the corresponding properties for the  cases of type II and type III respectively.

\begin{proposition}\label{prop:bian1}
For the case of type II,
the bi-anticanonical system of $S$ satisfies the following:
(i) the fixed component is $ C_1 + C_2 + \ol C_1 + \ol C_2$,
(ii) if we remove this fixed component,
the resulting system is base point free and $2$-dimensional,
(iii) if $\phi:S\to\CP^2$ is the induced morphism,
$\phi$ is of degree two, and the branch divisor is a quartic curve
which has two ordinary nodes,
(iv) the morphism $\phi$ maps the connected curves
$C_3\cup\ol C_1$ and $\ol C_3\cup C_1$
to the two nodes,
(v) $\phi$ maps both of the curves $C_2$ and $\ol C_2$ to the line
$l$ connecting the two nodes,
(vi) $\phi^{-1}(l)=C$.
\end{proposition}


%
\begin{proposition}\label{prop:bian2}
For the case of type III, the bi-anticanonical system of $S$ satisfies the following:
(i) the fixed component is $ C_1 + C_2 + C_3 + \ol C_1 + \ol C_2 + \ol C_3$,
(ii) after removing this,
the resulting system is base point free and $2$-dimensional,
(iii) the induced morphism $\phi:S\to\CP^2$  is of degree two, and the branch divisor is a quartic curve
with two cusps,
(iv)  $\phi$ maps the connected curves
$\ol C_4 \cup C_1\cup C_2$ and $ C_4 \cup \ol C_1\cup \ol C_2$
to the two cusps,
(v) $\phi$ maps both of the curves $C_3$ and $\ol C_3$ to the line
$l$ connecting the two cusps,
(vi) $\phi^{-1}(l)=C$.
\end{proposition}

\begin{figure}
\includegraphics{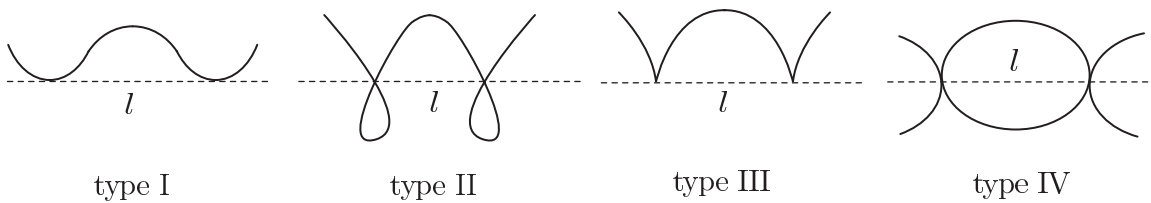}
\caption{The branch quartic curves, in relation
with the special line $l$.}
\label{fig:quartics}
\end{figure}

We omit proofs of these two propositions but
instead illustrate the branch quartic curve and the line $l$ as in Figure \ref{fig:quartics}.
(In the figure the branch quartic curves in the cases of 
type I and type IV are included from coherency of these
four types of the spaces.)

Going back to the twistor spaces, as consequences of these two propositions, we have

\begin{proposition}\label{prop:baseanti}
(i) For the case of type II, we have $\Bs\,|2F| = C_1
\cup C_2 \cup \ol C_1\cup \ol C_2$.
(ii) For the case of type III, we have $\Bs\,|2F| = C_1
\cup C_2 \cup C_3 \cup \ol C_1\cup \ol C_2\cup \ol C_3$.
\end{proposition}

\proof
These immediately follow from (i) of Propositions \ref{prop:bian1} and \ref{prop:bian2}, and surjectivity of the 
restriction map $H^0(Z, 2F)\to H^0(S, 2K_S^{-1})$
which was proved in \cite[Proposition 2.10]{HonDS4_3}.
(Note that the last surjectivity may be shown very readily
since we are on $4\CP^2$.)
\proofend

\vsp
We end this section by
the following property about reducible members of the pencil $|F|$,
which seems to have been well-understood (see Kreussler \cite{Kr98}).

\begin{proposition}
\label{prop:redF}
Let $Z$ be a twistor space on $n\CP^2$ satisfying $\dim |F| = 1$
and suppose that $\Bs\, |F|$ is a cycle  of rational curves,
written still by $C$.
Then the number of reducible members of the pencil $|F|$ is equal to
the half of the number of the components of $C$.
Further, all the reducible members are of the form
$S_i^+ + S_i^-$, where $S_i^{\pm}$ are degree-one divisors satisfying
$\ol S_i^+ = S_i^-$.
Furthermore, each of these members
divides the cycle $C$ into `halves'
in the sense that both $S_i^+\cap C$ and $S_i^-\cap C$ are connected
and  have equal components. 
\end{proposition}

In particular, $Z$ of type II (resp.\,III) has exactly three
(resp.\,four) reducible members of the pencil $|F|$.
We fix the indices of the reducible members by the property that
 $S_i^+\cap S_i^-$ (which is always a twistor line) goes through the point $C_i\cap C_{i+1}$,
 where we read $C_{i+1} = \ol C_1$ when $C_{i+1}$ does not exist.
 Also we make distinction between $S_i^+$ and $S_i^-$ by the property that 
 $S_i^+\supset \ol C_1$.
 These degree one divisors will play some role in our analysis
of the twistor spaces.

\section{Existence of some reducible anticanonical divisors}
In this section we prove that on the twistor spaces 
under consideration 
 there exist  some {\em irreducible} degree two divisors,
which do not belong to the fundamental system $|F|$.
Similarly to the case of type I studied in \cite{HonDS4_1},
these divisors will be a key for obtaining the defining equation of the branch divisor
of the anticanonical map.
However  in the present case 
for finding these divisors
we need different computations from the case
of type I given in \cite{HonDS4_1}.

Let $S^2H^0(F)$ be the subspace in $H^0(2F)$ generated by 
all elements in $H^0(F)$.
This is a 3-dimensional subspace of $H^0(2F)$.
\begin{proposition}
\label{prop:reducible}
Let $Z$ be  a Moishezon twistor space on $4\CP^2$
whose anticanonical map is (rationally) of  degree 2 over the image.
(i) 
If $Z$  is of type II,
then there exist distinct two anticanonical divisors $X_1+\ol X_1$ and $X_2 + \ol X_2$
on $Z$ which do not belong to  the subsystem $|S^2H^0(F)|$.
(ii) 
If $Z$ is of type III,
 there exists an anticanonical divisor $X+\ol X$ 
on $Z$ which does not belong to  the subsystem $|S^2H^0(F)|$.
\end{proposition}

We remark that $X_i+ \ol X_i\not\in |S^2H^0(F)|$ implies
irreducibility of $X_i$ and $\ol X_i$. 
We mention that if $Z$ is of type I (resp.\,type IV) there exist {\em three} 
(resp.\,{\em no}) such anticanonical divisors
as shown in \cite[Proposition 4.1]{HonDS4_1} (resp.\,\cite[Theorem 5.1, Eq.\,(100)]{HonDSn1}).
This difference for the number of such divisors will be reflected to a form of the
 defining equation of the branch divisor of the anticanonical map.

For the proof of Proposition  \ref{prop:reducible},
we need the following

\begin{lemma}\label{lemma:S}
Let $Z$ be as in Proposition  \ref{prop:reducible} and $S\in |F|$ a real irreducible member.
Then regardless of the type,
 there is a birational morphism $\epsilon : S\to \CP^1 \times \CP^1$
preserving the real structure,
such that the image  $\epsilon(C)$ is an anticanonical curve on $\CP^1\times\CP^1$,
and such that $\epsilon (C_1),\epsilon(\ol C_1)\in |\mathscr O(1,0)|$ and 
$\epsilon (C_2),\epsilon(\ol C_2)\in |\mathscr O(0,1)|$.
\end{lemma}

As a proof of this lemma is somewhat tedious to write down, we just mention that 
it suffices to
notice that if $A$ is any anticanonical curve on
a non-singular surface
and $D$ is a $(-1)$-curve satisfying $D\not\subset A$, then $A$ must intersect $D$ 
at a unique smooth point of $A$ and the intersection is transversal.

\vsp
\noindent
{\em Proof of Proposition \ref{prop:reducible} (i).}
Let $\epsilon: S\to \CP^1 \times \CP^1$ be the birational morphism as in Lemma 
\ref{lemma:S}.
This canonically determines a set of effective curves
$\{e_i,\ol e_i\set 1\le i\le 4\}$ 
satisfying $(e_i,e_j)_S= -\delta_{ij}$ and $(e_i,\ol e_j)_S=0$ 
for any $i$ and $j$.
(Since $\epsilon$ might involve a blowup at an infinitely near point,
$e_i$ and $\ol e_i$ can be reducible in general.)
From the self-intersection numbers of the components of the cycle $C$ in the case of type II, 
we can suppose that 
\begin{align}\label{int1}
e_3 = \ol C_3, \,\, (e_1,C_1) = (e_2,C_1) = 1, \,\, (e_4,C_2)=1.
\end{align}
Next let $t:Z\to 4\CP^2$ be the twistor fibration map, and let $\alpha_1,\alpha_2,\alpha_3,\alpha_4$
be orthonormal basis of $H^2(4\CP^2,\ZZ) $ which are uniquely determined by 
the property that $(t^*\alpha_i)|_S = e_i - \ol e_i$
in $H^2(S,\ZZ)$.

We first show that $h^0(F-t^*\alpha_1|_S) =1$.
In the following for simplicity we write $\alpha_i$ for the lift $t^*\alpha_i$,
and also write $(a,b):=\epsilon^*\mathscr O(a,b) \in H^2(S,\ZZ)$.
Then from \eqref{int1} and the bidegrees in Lemma \ref{lemma:S} we have the following relations in $H^2(S,\ZZ)$:
\begin{align*}
C_1 & = (1,0)- e_1-e_2-e_3,\\
C_2 & = (0,1) - \ol e_3 - e_4.
\end{align*}
By using this, we compute as
\begin{align}
(F - \alpha_1, C_1)_Z & = 
( K_S^{-1} - (e_1 - \ol e_1),\,   (1,0)- e_1-e_2-e_3 )_S\notag\\
& = ( (2,2) - 2e_1 - e_2 - e_3 - e_ 4 - \ol e_2 - \ol e_3 - \ol e_4 ,\,  (1,0)- e_1-e_2-e_3 )_S \notag \\
& = 2 + (-2 - 1 - 1) = -2.
\end{align}
Hence the curve $C_1$ is a fixed component of
the system $|(F-\alpha_1)|_S|$.
In a similar way we further find
$((F-\alpha_1)|_S - C_1,\, C_2)_S = -1,$
meaning $C_2$ is also a fixed component of the same system.
We then have
\begin{align}\label{subtract1}
(F-\alpha_1)|_S - C_1 - C_2 = (1,1) - e_1 - \ol e_2 - \ol e_4.
\end{align}
Now, since the cycle $C$ on $S$ consists of exactly 6 components,
for $i=1,2$, the points $\epsilon(e_i)$ and $\epsilon (\ol e_i)$ 
respectively belong to $\epsilon(C_1)$ and $\epsilon(\ol C_1)$,
{\em which are not  the singular points of the cycle $\epsilon(C)$.}
By the same reason, $\epsilon(e_4)$ and $\epsilon(\ol e_4)$ are  points on
$\epsilon(C_2)$ and $\epsilon(\ol C_2)$ respectively,
which are not the singular points of the cycle $\epsilon(C)$.
From these it follows that there exists a unique member of the 
linear system of \eqref{subtract1},
and that the member does not contain any of the irreducible
components of the cycle $C$.
(Such a member is exactly
the strict transform of the $(1,1)$-curve going through the three points
$\epsilon(e_1), \epsilon(\ol e_2)$ and $\epsilon( \ol e_4)$.)
Thus we get $h^0(F-\alpha_1|_S) =1$, as claimed.
In the same manner, we obtain 
$(F - \alpha_2|_S,\, C_1)_S<0$, $((F-\alpha_2)|_S - C_1,\,  C_2)_S <0$
and  $(F-\alpha_2)|_S - C_1 - C_2 = (1,1) - \ol e_1 - e_2 -\ol e_4$,
which again imply $h^0(F-\alpha_2|_S) =1$
and that the unique member of $|F-\alpha_2|_S|$ does not
contain any of the irreducible
components of the cycle $C$.

Next let $s\in H^0(F)$ be an element such that $(s)=S$,
and for $i=1,2$ we consider the obvious exact sequence 
\begin{align}\label{ses:4}
 0 \,\lra\,
 \mathscr O_Z(-\alpha_i) \,\stackrel{\otimes s}{\lra}\, F \otimes \mathscr O_Z(- \alpha_i) \,\lra\, ( F- \alpha_i) |_S \,\lra\, 0.
\end{align}
By Riemann-Roch formula we have $\chi(\mathscr O_Z(-\alpha_i))=0$ 
and by Hitchin vanishing 
\cite{Hi80} we have $H^2(\mathscr O_Z(-\alpha_i))=0$.
Also $H^0(\mathscr O_Z(-\alpha_i))=H^3(\mathscr O_Z(-\alpha_i))=0$ by trivial reason.
Hence we obtain $H^1 (\mathscr O_Z(-\alpha_i))=0$.
Thus by the cohomology exact sequence of \eqref{ses:4} we obtain $h^0(F-\alpha_i) = 
h^0( F-\alpha_i|_S)$.
Therefore we obtain $h^0(F-\alpha_i) = 1$ for $i=1,2$.

Let $x_i\in H^0(F-\alpha_i)$ be a non-zero element for $i=1,2$.
Then $\ol x_i:= \ol {\sigma^*x_i}$,
where $\sigma$ is the real structure of $Z$, is a non-zero element of $H^0 ( F+\alpha_i)$.
Hence the product $x_i\ol x_i$ belongs to $H^0(2F)$.
For finishing a proof we have to show 
$x_i\ol x_i\not\in S^2 H^0(F)$.
For this, we recall from the above argument that the divisor $(x_1|_S)$ contains the unique curve of the system 
\eqref{subtract1},
and that the curve does not contain components of 
the cycle $C$.
In particular the curve  $(x_1\ol x_1|_S)$ is not contained in $C$.
On the other hand, any $x\in S^2H^0(F)$ clearly satisfies
$x|_S=0$ or $(x|_S) = 2C$.
Thus we conclude $x_1\ol x_1\not\in S^2 H^0(F)$.
By the same argument we also obtain $x_2\ol x_2\not\in S^2 H^0(F)$.
Hence by letting $X_i=(x_i)$ and $\ol X_i = ( \ol x_i)$
for $i=1,2$, we finish a proof of Proposition \ref{prop:reducible} (i).
\proofend

\vspace{2mm}
The idea for (ii) being similar but again subtle, we give an outline:

\vsp
\noindent
{\em Proof of Proposition \ref{prop:reducible} (ii).}
Let $\epsilon : S\to \CP^1 \times \CP^1 $, $\{e_i,\ol e_i\set 1\le i\le 4\}$ and $\{\alpha_i\set 1\le i\le 4\}$ have the same meaning as in (i).
Then this time we can suppose
\begin{align}
C_1 = (1,0) - e_1 - e_2 - e_3, \,\,
C_2 = (0,1) - \ol e_3 - e_4, \,\,
C_3 = \ol e_3 - \ol e_2,\,\, C_4 = \ol e_2
\end{align}
in $H^2(S,\ZZ)$.
By computing intersection numbers we can show  that $|( F-\alpha_1)|_S|$ has $C_1+C_2+C_3$
as fixed components, and that 
$
(F - \alpha_1)|_S -C_1-C_2 - C_3 = (1,1) - e_1 - \ol e_3 - \ol e_4.
$ 
Further this system has a unique member, which is irreducible.
Then by the same argument for (i) we deduce $|F-\alpha_1|$ has a unique member 
$X$,
and that $X + \ol X\in |2F|$ and $X + \ol X\not\in |S^2H^0(F)|$.
\proofend


\section{Analysis of the anticanonical system of the twistor spaces}
\label{s:analysis}
In this section we analyze structure of the anticanonical map
of the twistor spaces.
For  type I twistor spaces, as showed in \cite[Proposition 3.2]{HonDS4_1},  
the base locus of the anticanonical
system $|2F|$  can be eliminated  by just blowing up 
the two $(-3)$-curves in the cycle $C$
(see \eqref{I}).
However, this is never the case for  type II and type III.
In this section we explicitly provide a full elimination of 
the base locus $\Bs\,|2F|$
for these two cases.
This is a core part of our analysis, and  indispensable for obtaining a constraint for the defining equation
of the branch divisor of the anticanonical map.
(In this section we do not need the results in the previous section.)

The elimination we take here is different from  \cite{HonDS4_1} for the case of type I,
in that we first blowup the whole of the cycle $C$, not the base curves themselves.
(So it is not a `minimal' elimination.) 
Though this yields a lot of  ordinary double points,
this provides $Z$ a structure of fibration, 
and this makes
much easier to  keep track of the base locus of the linear system,
which is otherwise quite difficult.

\subsection{The case of type II}
\label{ss:elimII}
Let $Z$ be a twistor space of $4\CP^2$ whose anticanonical map is degree two,
which is  of type II.
As before let
$\Phi:Z \to \CP^4$ be the anticanonical map,
$S\in |F|$ a real irreducible member,
and $C$ the unique anticanonical curve of $S$ (i.e.\,$C=\Bs \, |F|$).
$C$ is a cycle of six rational curves.
Let $f:Z\to \CP^1$ be the rational map associated with the pencil $|F|$.
The last $\CP^1$ can be naturally identified with the conic $\Lambda$
through the diagram \eqref{016}.
The map $f$ has the cycle $C$ as the indeterminacy locus.
Let $\mu_1: Z_1\to Z$ be the blowing-up at $C$, and
$E_i$ and $\ol E_i$ $(1\le i\le 3)$ the exceptional divisors
over $C_i$ and $\ol C_i$ respectively.
Then thanks to the fact that $C$ is a cycle,
all these exceptional divisors are isomorphic to
$\CP^1\times\CP^1$.
Write the composition 
$Z_1 \stackrel{\mu_1}{\to} Z \stackrel{f}{\to}\CP^1 $
by $f_1$.
Then since $F|_S\simeq \mathscr O_S(C)$,
 $f_1: Z_1 \to\CP^1 $ is 
a morphism.
By $\mu_1$ any  fiber of $f_1$ can be biholomorphically identified with a member of the pencil $|F|$.
Hence by Proposition \ref{prop:redF},
 $f_1$ has exactly three reducible fibers.
Let $\lambda_1,\lambda_2,\lambda_3\in \CP^1$ be the points
such that  $f^{-1}(\lambda_i) = S_i^+ \cup S_i^-$.
We put $L_i = S_i^+\cap S_i^-$ for $1\le i\le 3$.
These are twistor lines.

\begin{figure}
\includegraphics{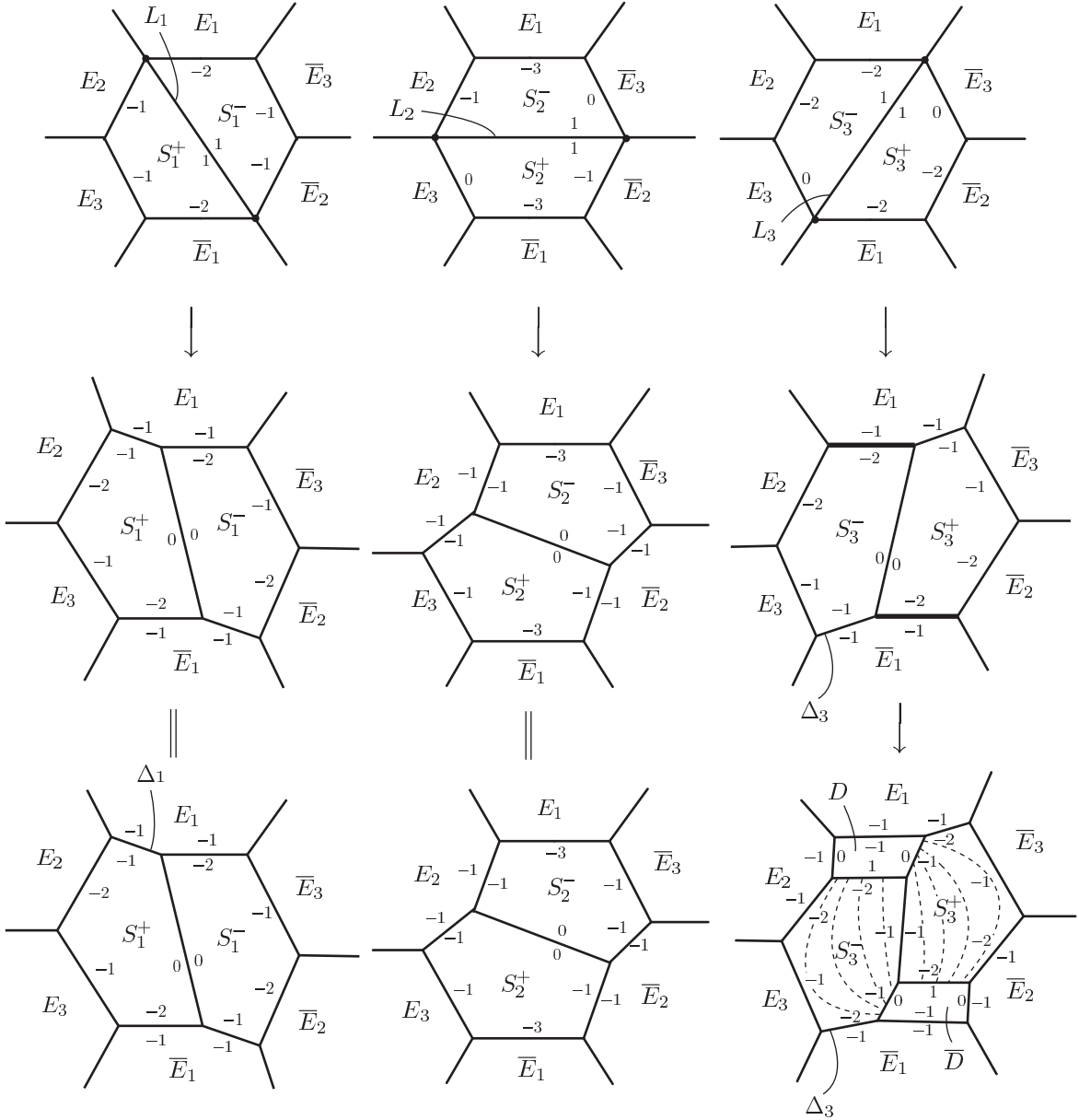}
\caption{The blowups for the case of type II.
The first, second, and the third rows are  
contained in $Z_1$,
$Z_2$ and $Z_3$ respectively.}
\label{fig:blowups_II}
\end{figure}

For simplicity we use the same letters to denote divisors and twistor lines in
the original space $Z$ and their strict transforms in $Z_1$.
In the first row of Figure \ref{fig:blowups_II} we illustrate $Z_1$ in a neighborhood of
each reducible fiber ($S_i^+\cup S_i^-$) of $f_1$.
As a computation using local coordinates shows,
 $Z_1$ has an ordinary node at the point
 where four faces meet, and these are indicated as dotted points.
(On each  reducible fiber there are two such points.)
We take  small resolutions for these six nodes as displayed in
the second row of Figure \ref{fig:blowups_II}.
This inserts two $\CP^1$-s in each reducible fiber of $f_1$,
and does not change any other part.
Let $ Z _2 $ be the resulting non-singular space,
and write $\mu_2:Z_2\to Z$ and $f_2:Z_2\to \CP^1$
for the compositions
$Z_2 \to Z_1 \stackrel{\mu_1}\to Z$ and $ Z_2\to Z_1
\stackrel{f_1}\to \CP^1$ respectively.
We again use the same letter for strict transforms into $Z_2$.
Then because of Proposition \ref{prop:bian1} (i), 
the pullback system $| \mu_2^*(2F) |$ has $E_1 + E_2 +  \ol E_1 + \ol E_2$ 
as fixed components at least.
Hence we put $$\mathscr L_2 : = \mu_2 ^* (2F) - (E_1 + E_2 +  \ol E_1 + \ol E_2).$$
Once we fix a section $v$ of
the line bundle
$\mathscr O_{Z_2}(E_1 + E_2 +  \ol E_1 + \ol E_2)$,
the birational morphism $\mu_2:Z_2\to Z$ naturally induces a  isomorphism
$H^0(2F)\simeq H^0(\mathscr L_2)$.
Moreover noticing the basic relation 
\begin{align}
f_2^*\mathscr O_{\CP^1} (1) \simeq \mu_2^*F - \sum_{1\le i\le 3} (E_i + \ol E_i) ,
\end{align}
we obtain 
\begin{align}\label{basic2}
\mathscr L_2 \simeq f_2 ^* \mathscr O(2) + E_1 + E_2 + 2 E_3 + \ol E_1 + \ol E_2 + 2 \ol E_3.
\end{align}
Therefore for any fiber $S_{\lambda} := f_2^{-1} (\lambda)$, 
we have 
$
\mathscr L_2 |_ {S_{\lambda}} \simeq  \mathscr O_{S_{\lambda}}( E_1 + E_2 + 2 E_3 + \ol E_1 + \ol E_2 + 2 \ol E_3).
$ 
This is useful for computing the base locus of $ |\mathscr L_2|$.
In particular if $S$ is a non-singular member of the pencil $|F|$,
 this isomorphism identifies $|\mathscr L_2|_{S}|$ with
the system
$|C_1 + C_2 + 2 C_3 + \ol C_1 + \ol C_2 + 2 \ol C_3|$,
which is  the movable part of $|2K_S^{-1}|$
by Proposition \ref{prop:bian1} (ii). 

\begin{lemma}\label{lemma:base1}
The linear system $|\mathscr L_2|$ on $Z_2$ has the following properties:
(i)  $|\mathscr L_2|$ has no  fixed component,
(ii) the two smooth rational curves 
$(S_3 ^- \cap E_1)$ and  $ (S_3 ^+ \cap  \ol E_1)$ are 
base curves of $|\mathscr L_2|$.
(In Figure \ref{fig:blowups_II}, these are displayed as bold segments.
Also, at this stage we do not prove these are all base points
of $|\mathscr L_2|$, although this is actually the case.)
%
\end{lemma}

\proof
For (i), by Proposition \ref{prop:bian1} (i), it is enough to show that the exceptional divisors $E_i$ and $\ol E_i$, $i=1,2$, are not a fixed component of $|\mathscr L_2|$.
For any smooth fiber $S$ of $f_2$ 
we have the following commutative diagram:
$$
\begin{CD}
H^0 (Z_2,\mathscr L_2) @>{\sim}>> H^0(Z,2F)\\
@VVV @VVV\\
H^0(S, \mathscr L_2|_S) @>{\sim}>> H^0(S, 2K^{-1}),
\end{CD}
$$
where the upper horizontal isomorphism is 
the composition of a multiplication
of the section $v$ and  the isomorphism
$H^0(Z_2, \mu_2^*(2F))\simeq H^0(Z, 2F)$,
and the lower horizontal arrow is its restriction to $S$, which is also isomorphic by \eqref{basic2} and Proposition 
\ref{prop:bian1}.
Further, 
the right vertical arrow is  surjective (\cite[Proposition 2.10]{HonDS4_3}).
Hence so is the left  arrow.
Therefore, since the linear system 
$|\mathscr L_2|_S|
\simeq |C_1 + C_2 + 2 C_3 + \ol C_1 + \ol C_2 + 2 \ol C_3|$
is base point free by Proposition \ref{prop:bian1} (ii),
we obtain that $|\mathscr L_2|$ has no fixed point on the smooth fiber $S$.
Hence any $E_i$ and $\ol E_i$ cannot be  a fixed component of $|\mathscr L_2|$.

For (ii), we temporary write $C_i:= S_3^- \cap E_i
\,(\subset Z_2)$ for $1\le i\le 3$,
and also let $\Delta_3 = S_3^-\cap \ol E_1$
(which is an exceptional curve of the small resolution $Z_2\to Z_1$).
Then by \eqref{basic2}, we compute
\begin{align*}
(\mathscr L_2 , {C_1} )_{Z_2}  &= (E_1 + E_2 + 2 E_3 + \ol E_1 + \ol E_2 + 2 \ol E_3,\, C_1)_{Z_2} \\
&= (C_1 + C_2 + 2C_3+\Delta_3,\, C_1)_{S_3^-} \\
&= -2 + 1 = -1.
\end{align*}
Hence $C_1$ is a base curve of $|\mathscr L_2|$.
By reality $\ol C_1 = S_3^+\cap \ol E_1$ is also a
base curve of $|\mathscr L_2|$.
\proofend

\vsp
By Lemma \ref{lemma:base1}, we let $\mu_3:Z_3\to Z_2$ 
be the blowup at the base curves $(S_3 ^- \cap E_1) \cup (S_3 ^+ \cap  \ol E_1)$
and $D$ and $\ol D$ the exceptional divisors respectively.
(See the lower right in Figure \ref{fig:blowups_II}.)
$D$ and $\ol D$ are isomorphic to $\mathbb P(\mathscr O(1)\oplus \mathscr O)$ over $\CP^1$.
We put $$\mathscr L_3 := \mu_3^*\mathscr L_2 - (D + \ol D).$$
The blowup $\mu_3$ induces an isomorphism $H^0(\mathscr L_3)\simeq H^0(\mathscr L_2)$.
Let $f_3:Z_3\stackrel{\mu_3}\to Z_2 \stackrel{f_2}\to\CP^1$ be 
the composition. Then
 the relation \eqref{basic2} is valid for $\mathscr L_3$ if we replace $f_2$ with $f_3$.
(But note that on $Z_3$ the surfaces $S_3^-$ and $ E_1$ are separated by $D$
as in Figure \ref{fig:blowups_II}.)
The next proposition implies that the blowup $\mu_3$ terminates an elimination
of  the base locus of the original system $|2F|$
on the twistor space.

\begin{proposition}\label{prop:free1}
The system $|\mathscr L_3|$ on $Z_3$ is base point free.
\end{proposition}

\proof
Again we use the same letter to mean the strict transforms into $Z_3$.
By computing the degrees of 
the line bundle $\mathscr L_3$ over
the curves on $E_1$ and $E_3$ seen in Figure \ref{fig:blowups_II},
it is possible to show
 that the restrictions of $\mathscr L_3$ 
to $E_1$ and $E_3$ are both {\em trivial}.
On the other hand,  the system $|\mathscr L_3|$ also does not
have any $E_i$ or $\ol E_i$ as a fixed component by
Lemma \ref{lemma:base1} (i).
These imply that $E_1$ and $E_3$ must be  disjoint from $\Bs\,|\mathscr L_3|$.
In particular, with the aid of Proposition \ref{prop:baseanti} (i),
we have  
$$\Bs\,|\mathscr L_3| \subset E_2\cup \ol E_2
\cup D\cup \ol D.$$
Next we see that $(\Bs\,|\mathscr L_3| ) \cap E_2 = \emptyset$.
Since $|\mathscr L_2|$ has no fixed point on any smooth fiber of $f_2$
as shown in the proof of Lemma \ref{lemma:base1}, 
the same is true for $|\mathscr L_3|$ for any smooth fiber
of $f_3$.
On the other hand, we can deduce that 
the system $|\mathscr L_3|_{E_2}|$ is a pencil without 
a base point,
whose members are sections of the natural projection
$E_2\to \Lambda$ (which is the restriction of $f_3:Z_3\to\CP^1$).
These imply that $|\mathscr L_3|$ does not have a base point on $E_2$ too.

Finally we see that 
$(\Bs\,|\mathscr L_3| ) \cap D = \emptyset$.
For this we first notice from \eqref{basic2} and 
Figure \ref{fig:blowups_II} that 
the restriction $\mathscr L_3|_D$ is isomorphic to 
a pullback of $\mathscr O(1)$ by the blowdown $D\to \CP^2$.
Therefore in order to show $(\Bs\,|\mathscr L_3| ) \cap D = \emptyset$, it suffices to prove that 
the rational map $\Phi_3$ associated to 
$|\mathscr L_3|$ satisfies $\dim\Phi_3(D)=2$.
We show this by proving that $\Phi_3(D)$ contains two distinct 
lines.

From the construction $\Phi_3$ factors as $Z_3\to Z\stackrel{\Phi}\to \CP^4$, where $Z_3\to Z$ is the composition of all blowups.
Therefore the image of $\Phi_3$ is the scroll $Y=\pi^{-1}({\Lambda})$.
Moreover 
each fiber of the fibration $f_3:Z_3\to\CP^1$ is mapped
(by $\Phi_3$) to a fiber plane of the scroll $Y\to \Lambda$
(see the diagram \eqref{016}).
Next we show that $\Phi_3(S_3^-)$ is  a line.
 By \eqref{basic2}, temporary writing $ C_i=S_3^-\cap  E_i$ for $1\le i\le 3$ on $Z_3$ and recalling that 
we have subtracted the divisor $D$ when defined $\mathscr L_3$, we have 
\begin{align}\label{fa89fw}
\mathscr L_3|_{S_3^-}\simeq \mathscr O_{S_3^-} (  C_2 + 2 C_3 + \Delta_3).
\end{align}
From this it is easy to see that
the system 
$|\mathscr L_3|_{S_3^-}|$ is a pencil without
a base point, and the associated morphism
$S_3^-\to \CP^1$ has the curves
$D\cap S_3^-$ and $\ol D\cap S_3^-$ as sections.
(In the lower right of Figure \ref{fig:blowups_II},
fibers of the morphism $S_3^-\to \CP^1$ are indicated 
as dotted curves.)
In particular, if
the restriction 
$H^0(Z_3,\mathscr L_3)\to H^0(S_3^-, \mathscr L_3|_{S_3^-})\simeq\CC^2$ is not surjective, 
the system $|\mathscr L_3|$ has
base points along a fiber of the morphism
$S_3^-\to \CP^1$.
But this cannot happen 
since we already know $\Bs\,|\mathscr L_3|\subset
D\cup \ol D$.
Therefore the above restriction is surjective.
This implies that $\Bs\,|\mathscr L_3|\cap
S_3^-=\emptyset$ and $\Phi_3(S_3^-)$ is a line.
Hence $\Phi_3(S_3^-\cap D)$ is also the same line.
In particular $\Phi_3(D)$ contains the line
$\Phi_3(S_3^-)$.
Also, this line can be written as 
$\Phi_3(S_3^-\cap \ol D)$
(since $S_3^-\cap \ol D$ was a section of the 
morphism $S_3^-\to\CP^1$).

By the real structure, $\Phi_3(S_3^+)$ is also a line,
and this can also be written as $\Phi_3(S_3^+\cap  D)$.
This implies that $\Phi_3(D)$ contains 
the line $\Phi_3(S_3^+)$ too.
If $\Phi_3(S_3^+) = \Phi_3(S_3^-)$,
the image $\Phi(S_3^-\cup S_3^+)$ would be a 1-dimensional linear subspace of the plane $\pi^{-1}_3(\lambda_3)$,
which cannot happen.
Thus the two lines  $\Phi_3(S_3^+)$ and $\Phi_3(S_3^-)$
are distinct.
Hence $\Phi_3(D)$ contains two lines.
Therefore $\Phi_3(D)=\CP^2$, and we finally obtain 
$(\Bs\,|\mathscr L_3| ) \cap D = \emptyset$.
Thus we conclude $\Bs\,|\mathscr L_3|=\emptyset$,
as claimed.  
\proofend

\vsp
For restrictions of $\Phi_3$ to fibers of the morphism $f_3:Z_3\to\CP^1$, we have the following.

\begin{lemma}\label{lemma:anti2}
(i) If $S=f_3^{-1}(\lambda)$ is a smooth fiber of $f_3$, the restriction $\Phi_3|_S$ is of degree two over the plane
$\pi^{-1}(\lambda)$.
(ii) On two reducible fibers $S_1^+\cup S_1^-$ and $S_2^+\cup S_2^-$,
$\Phi_3$ is birational over the plane
$\pi^{-1}(\lambda_1)$ and $\pi^{-1}(\lambda_2)$ 
respectively on each irreducible component.
Further, the image of the twistor line $L_i= S_i^+ \cap S_i^-$ is a conic in the plane.
(iii) $\Phi_3$ contracts each of $S_3^+$ and $ S_3^-$ to a line in the plane.
Further, these two lines are distinct.
\end{lemma}

In terms of the original map $\Phi:Z\to\CP^4$, the above (iii) means that 
the anticanonical map contracts $S_3^+$ and $S_3^-$ to lines.
This is a major difference between the case of type I, 
where the anticanonical map does not contract any divisor \cite[Proposition 3.6]{HonDS4_1}.

\vsp\noindent
{\em Proof of Lemma \ref{lemma:anti2}.}
(i) is obvious because on such $S$ the restriction $\Phi_3|_S$ can be identified with 
the bi-anticanonical map of $S$, which is degree two over a plane by Proposition \ref{prop:bian1} (iii).
For (ii) in the case $i=1$, on $Z_3$ 
we temporary put $C_2=S_1^+\cap E_2, C_3 = S_1 ^+ \cap E_3$ and $\ol C_1= S_1^+\cap \ol E_1$.
We also put $\Delta_1=S_1^+\cap E_1$
(see lower left in Figure \ref{fig:blowups_II}).
Then by the isomorphism \eqref{basic2} we have $\mathscr L_3|_{S_1^+} \simeq \mathscr O_{S_1^+}
(\Delta_1 + C_2 + 2 C_3 + \ol C_1).$
By standard computations it is possible to show that:
\begin{gather*}
h^0 (S_1^+,\mathscr O
(\Delta_1 + C_2 + 2 C_3 + \ol C_1 )) = 3,\\
(\Delta_1 + C_2 + 2 C_3 + \ol C_1 )^2 = 1 {\text{ on $S_1^+$}},\\
\Bs\,|\Delta_1 + C_2 + 2 C_3 + \ol C_1| =\emptyset,
\end{gather*} and 
also the induced morphism $S_1^+\to \CP^2$ is birational.
The restriction $\Phi_3|_{S_1^+}$ is nothing but the rational map 
associated to the image of the restriction map  
$H^0(Z_3, \mathscr L_3) \to 
H^0(S_1^+, \mathscr O_{S_1^+}
(\Delta_1 + C_2 + 2 C_3 + \ol C_1 ))$.
The last image cannot be $0$ or $1$-dimensional since $\Bs\,|\mathscr L_3| \neq \emptyset$.
Also, it cannot be $2$-dimensional since 
in that case $|\mathscr L_3|_{S_1^+}|$ would have a base point because
 $(\mathscr L_3|_{S_1^+})^2 = 1$ as  above.
 Thus we conclude that $\Phi_3|_{S_1^+}$ is birational.
By  similar computations for which we omit, we can also show that 
$\Phi_3|_{S_2^+}$ is birational over the plane
$\pi^{-1}(\lambda_2)$.

(iii) is already shown in the final part of the  proof of Proposition \ref{prop:free1}.
\proofend 

\vsp
The morphism $\Phi_3$ maps
the exceptional divisors of the blowups as follows.
\begin{proposition}\label{prop:gfshu}
Let $\Phi_3:Z_3\to \CP^4$ be the morphism associated to $|\mathscr L_3|$ as before.
Then
(i) by $\Phi_3$, the divisors $E_1$ and $\ol E_3$ are mapped to one and the  same point 
on the singular line of the scroll $Y$.
The same is true for $\ol E_1$ and $ E_3$.
(ii) by $\Phi_3$, $E_2$ and $\ol E_2$ are mapped onto the singular line of $Y$.
(iii) by $\Phi_3$, $D$ and $\ol D$ are mapped birationally to the plane $\pi^{-1}(\lambda_3)$.
\end{proposition}

\proof
We first show (ii).
As is remarked in the proof of Proposition 
\ref{prop:free1}, the restriction $|\mathscr L_3 |_{E_2}|$ is a  pencil without a base point.
Therefore by the same reason for $S_3^+$ in Lemma \ref{lemma:anti2}, 
$E_2$ is mapped to a line by $\Phi_3$.
We show this line must be the singular line of the scroll. For any non-singular fiber $S$ of $f_3$, $\Phi_3|_S$ is naturally identified with the 
bi-anticanonical map of $S$, and $S\cap E_2$ is exactly the component $C_2$, 
which is mapped to a line by Proposition \ref{prop:bian1} (v).
Since this line is exactly $\Phi(C_2)$,
 this is independent of the choice of $S$
Therefore
the line must be the intersection of fiber planes
of the projection $\pi:\CP^4\to\CP^2$.
Thus the $\Phi_3(E_2)$ must be the singular line of $Y$.
Hence since the singular line is real, $\ol E_2$ 
is also mapped to the same line, and we get (ii).

Since $\mathscr L_3$ is trivial on $E_1$ and $\ol E_3$,
each of these are mapped to a point by $\Phi_3$.
Since $E_1\cap \ol E_3\neq\emptyset$, these points
must coincide.
Further, as $E_1\cap E_2\neq\emptyset$ 
the last point must belong to the singular line of the scroll.
By real structure we obtain the same conclusion for
$\ol E_1$ and $E_3$, and we obtain (i).
Finally (iii) is already shown in the course of the proof
of Proposition \ref{prop:free1}.
\proofend

\subsection{The case of type III}
\label{ss:elimIII}
The elimination of the base locus of the anticanonical map 
in the case of type III can be done along  the same line as in the case of type II
(though more complicated).
So the description below is partially sketchy.

\begin{figure}
\includegraphics{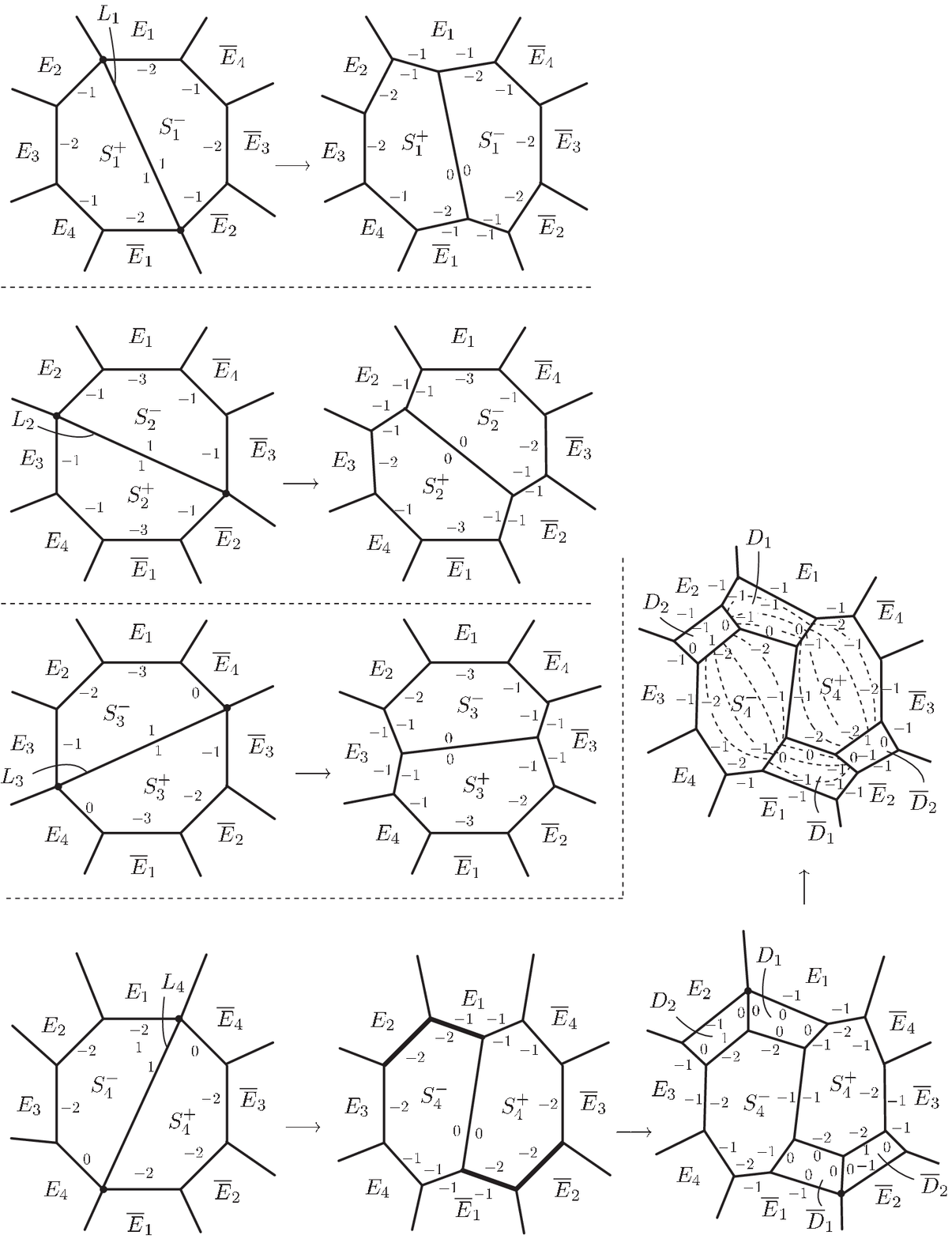}
\caption{The blowups for the case of type III.}
\label{fig:blowups_III}
\end{figure}

Let $Z$ be a twistor space on $4\CP^2$ with degree two anticanonical map, which is of type III.
Let $\mu_1:Z_1\to Z$ be the blowup at $C$, and $E_i$ and $\ol E_i$ $(1\le i\le 4$ this time)  the exceptional divisors over
$C_i$ and $\ol C_i$ respectively.
Then again $Z_1$ admits a morphism $f_1:Z_1\to \CP^1$
induced by the pencil 
$$\bigg|\mu_1^* F-\sum_{1\le i\le 4} (E_i + 
\ol E_i)\bigg|.
$$
This fibration $f_1$ has exactly four reducible fibers and 
they can be described as in the first column of
Figure \ref{fig:blowups_III}
in a neighborhood of the fibers.
Let $\lambda_i$, $1\le i\le 4$, 
be
the points on $\Lambda$ corresponding to the reducible fiber
$S_i^+\cup S_i^-$.
As in the case of type II, $Z_1$ has
exactly two ordinary nodes on each reducible fiber
(again indicated as dotted points).
For each of these eight nodes of $Z_1$ 
we take a small resolution as displayed in the figure.
Let $\mu_2:Z_2\to Z_1\stackrel{\mu_1}{\to} Z$ be the composition,
and put $$\mathscr L_2:=\mu_2^*(2F) - (E_1 + E_2 + E_3 + \ol E_1 
+ \ol E_2 + \ol E_3).$$
(See Proposition \ref{prop:bian2} (i).) 
Let $f_2:Z_2\to Z_1\stackrel{f_1}{\to}\CP^1$ be the composition.
We have a natural isomorphism
$H^0(2F)\simeq H^0(\mathscr L_2)$ and also, similarly to \eqref{basic2} in the case of type II, an isomorphism
\begin{align}\label{basic4}
\mathscr L_2 \simeq f_2 ^* \mathscr O_{\CP^1}(2) + E_1 + E_2 + E_3 + 2 E_4 + \ol E_1 + \ol E_2 + \ol E_3 + 2 \ol E_4.
\end{align}
Then analogously to Lemma \ref{lemma:base1} we have the following

\begin{lemma}\label{lemma:base2}
On $Z_2$, we have the following:
(i) the system $|\mathscr L_2|$ does not have a fixed component,
(ii) the following four smooth rational curves 
\begin{align}\label{4curves}
S_4 ^+ \cap \ol E_1,  \quad S_4 ^+ \cap \ol E_2,  \quad
S_4 ^- \cap  E_1,  \quad S_4 ^- \cap  E_2
\end{align} are 
base curves of $|\mathscr L_2|$.
(Note that the first two curves intersect and the same for the last two curves; see Figure \ref{fig:blowups_III}.)
%
\end{lemma}

\proof
(i) is completely analogous to the proof of Lemma \ref{lemma:base1} (i).
(We use Proposition \ref{prop:bian2} (ii) instead of Proposition \ref{prop:bian1} (ii).)
For (ii), from \eqref{basic4} and Figure \ref{fig:blowups_III}, the  intersection numbers of $\mathscr L_2$ with
the four curves \eqref{4curves} are respectively 
computed to be $-1,0,-1,0$.
These imply the claim (ii).
\proofend

Let $\mu_3:Z_3\to Z_2$ 
be the blowup at the base curves \eqref{4curves}
and $D_1$ and $D_2$ the exceptional divisors
over $S_4 ^- \cap  E_1$ and $ S_4 ^- \cap  E_2 $ respectively.
A careful computation shows 
\begin{align}\label{Ds}
D_1\simeq\CP^1\times\CP^1,
\quad
D_2\simeq\mathbb P(\mathscr O(1)\oplus\mathscr O).
\end{align}
$Z_3$ has ordinary nodes over the two singular points of the
center of $\mu_3$.
Let $Z_4\to Z_3$ be the small resolution of these two nodes
as displayed in Figure \ref{fig:blowups_III},
and $\mu_4:Z_4\to Z_3\to Z_2$ the composition.
Then by the last small resolution $D_1$ is blownup at a point, 
while $D_2$ remains unchanged
(see Figure \ref{fig:blowups_III}).
We again use the same letters to denote the strict transforms of these divisors into $Z_4$,
and put
$$\mathscr L_4:=\mu_4^*\mathscr L_2 - (D_1+ D_2 + \ol D_1 + \ol D_2).
$$
Then similarly to Proposition \ref{prop:free1}, we have
\begin{proposition}\label{prop:free2}
The system $|\mathscr L_4|$ on $Z_4$ is base point free.
\end{proposition}

\proof
By Lemma \ref{lemma:base2} (i) the divisors $E_i$ and $\ol E_i$ are not a fixed component of $|\mathscr L_4|$ for any $1\le i\le 4$.
Further by using the relation \eqref{basic4} and 
chasing Figure \ref{fig:blowups_III} it is possible 
to show that the line bundle $\mathscr L_4$ is 
trivial on $E_1, E_2$ and $E_4$.
These mean that $|\mathscr L_4|$ does not have 
a base point on $E_1, E_2$ and $E_4$.
Hence, with the aid of Proposition \ref{prop:baseanti} (ii),
we have
$$
\Bs\,|\mathscr L_4| \subset E_3\cup D_1\cup D_2 \cup 
\ol E_3\cup \ol D_1 \cup \ol D_2.
$$
Moreover, the system $|\mathscr L_4|_{E_3}|$ satisfies 
$\dim |\mathscr L_4|_{E_3}|=1$ and
the intersection number with a fiber class of 
the natural projection $E_3\to \Lambda$ is one.
Therefore, since $|\mathscr L_4|$ does not have 
a base point on a smooth fiber of the composition 
$Z_4 \to Z_1\to \Lambda$, we obtain that 
$|\mathscr L_4|$ does not have a base point on $E_3$ too.

Let $\Phi_4:Z_4\to Y\subset\CP^4$ be the map
associated to $|\mathscr L_4|$.
It remains to see that there is no base point on $D_1$ nor $D_2$.
By \eqref{Ds}, $D_2$ is biholomorphic to one point blownup
of $\CP^2$.
Further, again from Figure \ref{fig:blowups_III} and 
\eqref{basic4}, 
we can deduce that $\mathscr L_4|_{D_2}$ is isomorphic to the pullback
of $\mathscr O(1)$ by a blowdown $D_2\to \CP^2$.
Also it is possible to deduce that 
 the restrictions of $\mathscr L_4$
to the divisors $S_4^-$ and $D_1$  are
respectively isomorphic to the pullback of 
$\mathscr O_{\CP^1}(1)$ by  surjective morphisms
$S_4^-\to \CP^1$ and $D_1\to \CP^1$,
the curves $S_4^-\cap D_2$ and $S_4^-\cap\ol D_1$ are
sections of the morphism $S_4^-\to \CP^1$,
and that both $D_1\cap S_4^+$ 
and $D_1\cap D_2$ are sections of
the morphism $D_1\to\CP^1$.
(In Figure \ref{fig:blowups_III} fibers of these morphisms are
indicated by dotted curves.)
From these we obtain that the restrictions of
$\Phi_4$ to $S_4^-$ and ${D_1}$ are exactly 
the above morphisms to $\CP^1$, and that
the image $\Phi_4(S_4^-)$ and $\Phi_4(D_1)$ are lines
in the plane $\pi^{-1}(\lambda_4)$.
In particular we have $(\Bs\,|\mathscr L_4|) \cap D_1 
=\emptyset$.
We also have $\Phi_4(S_4^-) = \Phi_4(S_4^-\cap \ol D_1)$ and
$\Phi_4(D_1) 
= \Phi_4(D_1\cap D_2) = \Phi_4(D_1\cap S_4^+)$.
Hence by the real structure we have $\Phi_4(D_1\cap D_2)= \Phi_4(S_4^+)$.
Therefore the above two lines in $\pi^{-1}(\lambda_4)$
can be rewritten as $\Phi_4(S_4^-)$ and $\Phi_4(S_4^+)$,
which means that these are distinct.
Thus the image $\Phi_4(D_2)$ contains 
two distinct lines and hence 
$\Phi_4|_{D_2}:D_2\to \Phi_4(D_2)$ is just the blowdown
(of the $(-1)$-curve $D_2\cap E_2$).
In particular, $(\Bs\,|\mathscr L_4|)\cap D_2 = \emptyset$.
Hence we obtain $\Bs\,|\mathscr L_4| = \emptyset$.
\proofend

\vsp
By the proposition, the map $\Phi_4:Z_4 \to Y\subset\CP^4$
is a morphism.
For the images of  divisors on $Z_4$ by $\Phi_4$, we have
the following 

\begin{proposition}\label{prop:dsg84ga}
(i) For $i=1,2,3$, the restrictions of $\Phi_4$ to the  divisors
$S_i^{\pm}$ are birational over the plane $\pi^{-1}(\lambda_i)$,
(ii) $\Phi_4$ contracts $S_4^+$ and $S_4^-$ to lines in $\pi^{-1}(\lambda_4)$ which are mutually distinct,
(iii) for $i=1,2,4$, $\Phi_4$ contracts the divisors $E_i$
and $\ol E_i$ to points on the singular line $l$ of the scroll,
(iv) $\Phi_4$ maps $E_3$ and $\ol E_3$ to the line $l$,
(v) $\Phi_4(D_1) = \Phi_4(S_4^+)$ and 
$\Phi_4(\ol D_1) = \Phi_4(S_4^-)$,
(vi) the restrictions of $\Phi_4$ to $D_2$ and $\ol D_2$ are
birational over the plane $\pi^{-1}(\lambda_4)$.
\end{proposition}

\proof
The claims (ii), (iii), (v) and (vi) are already proved in 
the proof of Proposition \ref{prop:free2}.
(i) and (iv) can be shown in a similar way to 
Lemma \ref{lemma:anti2} (i) and Proposition \ref{prop:gfshu} (ii) respectively.
\proofend 
\subsection{The branch divisor of the anticanonical map}
Let $Z$ be a twistor space on $4\CP^2$, which is of type II or type III, and $\Phi:Z\to Y\subset\CP^4$
the anticanonical map as before.
We  write $\nu:\tilde Y\to Y$ for the blowup of $Y$
along the singular line $l$.
Since $\Lambda$ is a conic,
$\tilde Y$ is biholomorphic to 
the total space of the $\CP^2$-bundle $\mathbb P
(\mathscr O(2)^{\oplus 2}\oplus \mathscr O)\to\CP^1
\simeq\Lambda$.
We say that a morphism to $Y$ can be lifted to $\tilde Y$ 
if it factors through the resolution $\tilde Y\to Y$.
\begin{lemma}\label{lemma:lift}
(i) If $Z$ is of type II, the morphism $\Phi_3:Z_3\to Y$
obtained in Section \ref{ss:elimII} can be lifted to a morphism
$\tilde{\Phi}_3:Z_3\to\tilde Y$.
(ii) If $Z$ is of type III, the morphism ${\Phi}_4:Z_4\to Y$
obtained in Section \ref{ss:elimIII} can be lifted to a morphism
$\tilde{\Phi}_4:Z_4\to\tilde Y$.
\end{lemma}

\proof
Since this can be shown as in the proof of the corresponding
statement \cite[Proposition 3.4]{HonDS4_1},
we just give a sketch of the proof.
The indeterminacy locus of  the rational map $f:Z\to \CP^1$
associated to the pencil $|F|$ was eliminated  
by just blowing-up the cycle $C$,
and we always have a morphism to $\CP^1$ from the series of 
blownup spaces.
Via the natural identification between 
fibers of $\tilde Y\to \Lambda$ and planes in the scroll $Y$,
this gives the desired lifts of the morphisms $\Phi_3:Z_3\to Y$
and $\Phi_4:Z_4\to Y$.
\proofend

\vsp
As the original anticanonical map $\Phi:Z\to Y$ is (rational but) of 
degree two, the lifts are also of degree two.
For the branch divisor of these lifts, we have

\begin{proposition}\label{prop:branch}
Let $\tilde{\Phi}_3:Z_3\to \tilde Y$ and $\tilde{\Phi}_4:Z_4\to \tilde Y$
be the lifts as in Lemma \ref{lemma:lift}.
Then the branch divisor of these degree-two morphisms
are pull-back by $\nu:\tilde Y\to Y$ of a divisor belonging to the system $|\mathscr O_Y(4)|$, where $\mathscr O_Y(4):=\mathscr O_{\CP^4}(4)|_Y$.
\end{proposition}

Because this can be proved in a similar way to 
the case of type I (\cite[Proposition 3.4]{HonDS4_1}),
we just mention that it is enough to use
Propositions \ref{prop:bian1} and \ref{prop:bian2}
instead of  \cite[Propositions 2.2 and 2.3]{HonDS4_1}.

By Proposition \ref{prop:branch}, the branch divisors of the
morphisms $\Phi_3:Z_3\to Y$ and $\Phi_4:Z_4\to Y$ 
are cuts of $Y$ by  quartic hypersurfaces.
Needless to say, this hyperquartic is the most 
significant data for determining the twistor spaces.
In the next section we examine
defining equation of these hyperquartics.

\section{Defining equations of the branch  divisors}
As in the previous section let $Z$ be a Moishezon twistor space on $4\CP^2$
whose anticanonical map is degree two which is of type II 
or type III, and $\Phi_3:Z_3\to Y$ or $\Phi_4:Z_4\to Y$
respectively be the degree two morphisms which are obtained as 
a consequence of the elimination 
of the original anticanonical map $\Phi:Z\to Y$,
given in the last section.
We denote by $B\,(\subset Y)$ the branch divisor
of the morphisms.
By Proposition \ref{prop:branch} we know that $B$ is of the form
$Y\cap\mathscr B$, where 
$\mathscr B$ is a quartic hypersurface in $\CP^4$.

\subsection{Double curves on the branch divisor}
\label{ss:dc}
Let $H$ be a hyperplane in $\CP^4$.
The intersection $Y\cap H$ is either a plane or a cone over the conic $\Lambda$.
We say that a curve $\mathscr C$ on the branch divisor $B$ is a {\em double curve
with respect to $H$} if the hyperplane $H$ touches $B$ along the curve $\mathscr C$;
or more precisely if $B\cap H$ is a non-reduced curve
on a reduced surface $Y\cap H$ (i.e.\,a plane or a cone).
In this subsection we  find five double curves on $B$.

First we consider the case of type II.
For $i=1,2$ let $L_i=S_i^+\cap S_i^-$ be the 
intersection twistor line.
Then by Lemma \ref{lemma:anti2} the image $\Phi_3(L_i)$, which clearly coincides with $\Phi(L_i)$,  is a plane conic. We put $\mathscr C_i:= \Phi(L_i)$.
Then as $\Phi$ is birational on each of the components
$S_1^{\pm}$ and $S_2^{\pm}$ by the lemma,
$\mathscr C_1$ and $\mathscr C_2$ are double curves 
with respect to  hyperplanes containing 
the planes $\pi^{-1}(\lambda_1)$ and $\pi^{-1}(\lambda_2)$
respectively.
We call these two curves as {\em double  conics}.

On the other hand the situation is very different for
$S_3^+$ and $S_3^-$ by Lemma \ref{lemma:anti2} (iii) 
in that these are contracted to lines by $\Phi$.
We define 
$$\mathscr C_3:= \Phi(S_3^+\cup S_3^-).$$ 
By the lemma this is a union of two distinct lines in the plane
$\pi^{-1}(\lambda_3)$.
From the definition the curve $\mathscr C_3$ must also 
be a double curve with respect to  hyperplanes
containing $\pi^{-1}(\lambda_3)$.
We call $\mathscr C_3$  {\em a splitting double conic}.
This kind of a double curve did not appear in the case of 
type I and 
 will be significant in obtaining 
a strong constraint for defining equation of the branch divisor.
These three double conics $\mathscr C_1,\mathscr C_2$ and $\mathscr C_3$ 
are placed in $\CP^4$ in a way that the intersection
$\mathscr C_i\cap l$ (where $l$ is the singular line of the scroll
as before) consists of two points which are independent of $i$.
More invariantly, recalling that 
the bi-anticanonical map $\phi$ on $S$ contracts
the connected curves $\ol C_3\cup C_1$ and $ C_3\cup \ol C_1$
(Proposition \ref{prop:bian1} (iv)),
these two points on $l$ are the images 
$\Phi(\ol C_3\cup C_1)$
and $\Phi(C_3\cup \ol C_1)$ respectively.

Next recalling that the anticanonical system $|2F|$ has
two distinguished members $X_1 + \ol X_1$
and $X_2 + \ol X_2$ by Proposition \ref{prop:reducible}
(i) we define $H_i$ ($i=4,5$) to be the unique hyperplane in $\CP^4$
which satisfies $\Phi^{-1}(H_i) = X_{i-3} + \ol X_{i-3}$.
We now show that $X_1$ and $X_2$ are not contracted to a curve
or a point by $\Phi$.
For this as in the proof of Proposition \ref{prop:reducible}
for a generic member $S\in |F|$ the intersection
$S\cap X$ contains a curve which does not contain any component
of the cycle $C$.
By using $K_S^2=0$, it is easy to show that 
the curves $C_1,C_3, \ol C_1$  and $\ol C_3$ are
all curves that are contracted to a point by the 
bi-anticanonical map $\phi$ of $S$.
Hence recalling $\Phi|_S=\phi$, the above curve $X\cap S$ is not contracted to 
a point by $\Phi$.
Because $S$ is generic,
this implies that $X$ itself is also not contracted to 
a curve by $\Phi$, as claimed.
Therefore $\Phi|_{X_i}$ is birational, and 
$\Phi(X_i\cap \ol X_i)$ $(i=1,2)$ must be a 
double curve with respect to $H_{i+3}$.
We write $\mathscr C_4$ and $\mathscr C_5$
 for these double curves.
As $Y$ is quadratic and $B\in |\mathscr O_Y(4)|$, these
curves must be of degree four in $\CP^4$.
We call these {\em double quartic curves}.
These two curves intersect at four points on the plane
$H_4\cap H_5$.
So for the case of type II we have obtained 
three double conics (one of which is a splitting one)
and two double quartic curves.

For the case of type III, in a similar way, the three curves
$\mathscr C_i:= \Phi(L_i)$ (where $L_i=S_i^+\cap S_i^-$
and $1\le i\le 3$ this time)
are double conics with respect to  hyperplanes containing the plane
$\pi^{-1}(\lambda_i)$.
Also if we put $\mathscr C_4:= \Phi(S_4^+\cup S_4^-)$,
 this becomes a splitting double conic
in the above sense by Proposition \ref{prop:dsg84ga} (ii).
These four double conics intersect the singular line $l$
at  two points which are independent of the choice of the conic.
These two points are nothing but the image 
under $\Phi$ of the connected curves
$\ol C_4\cup C_1\cup C_2$ and $ C_4\cup \ol C_1\cup \ol C_2$.
Further recalling Proposition \ref{prop:reducible} (ii) and letting $H_5$  be 
the hyperplane satisfying $\Phi^{-1}(H_5)= X+ \ol X$,
$\mathscr C_5:=\Phi(X\cap \ol X)$ becomes 
a double quartic curve with respect to $H_5$
by the same reason for $\mathscr C_4$ and $\mathscr C_5$
in the case of type II.
Thus for the case of type III we have obtained 
four double conics (one of which is a splitting one)
and one double quartic curve.

We emphasize that for both of types II and III
double conics are lying on fiber planes of 
the scroll $Y\to\Lambda$, whereas  double quartic curves are lying on  cones which are hyperplane sections of the scroll.
This difference will be significant when we determine
 defining equation of the branch divisor of
the anticanonical map of the twistor spaces.

\begin{remark}
{\em We now display the number of double curves, including the
cases of type I and type IV.
For the case of type I, these are proved in \cite[Section 4.1]{HonDS4_1}. The case of type IV can be proved in a similar way.
(But note that in this case we did not give a full elimination of the base locus in \cite{HonDSn1}.)
\begin{center}
\begin{tabular}{|c||c|c|c|c|}
\hline
                & type I & type II & type III & type IV\\ \hline
double conics (splitting one)  & 2  (0) &    3 (1)&  4  (1)  & 5  (1)  \\ 
double quartic curves & 3      &    2    &   1      & 0\\ \hline
total number & 5     &    5   &   5      & 5\\ \hline
\end{tabular} 
\end{center}
Note that the number of the double conics is always equal to the 
number of the reducible members of the pencil $|F|$.
}
\end{remark}

\subsection{Quadratic hypersurfaces containing 
double curves}
Next  we investigate  hyperquadrics in $\CP^4$ containing
all these five double curves.

\begin{proposition}\label{prop:Q}
As before let $Z$ be of type II or type III,
$B\subset Y$ the branch divisor of the anticanonical map,
and $\mathscr C_i, \,1\le i\le 5$, the double curves on $B$.
Then there exists a hyperquadric in $\CP^4$ 
 which contains all these double curves and which is different
from the scroll $Y$.
\end{proposition}

\proof
For the case of type II, 
we consider quadratic hypersurfaces in $\CP^4$ which 
go through all the following 12 points:
\begin{enumerate}
\item[(a)] $\mathscr C_1\cap\mathscr C_2$ (consisting of 2 points),
\item[(b)] $\mathscr C_4\cap \mathscr C_5$ 
(consisting of 4 points),
\item[(c)]
$\mathscr C_1\cap \mathscr C_4$,
and one of the two points $\mathscr C_1\cap \mathscr C_5$
(consisting of 3 points),
\item[(d)] $\mathscr C_2\cap \mathscr C_4$
(consisting of 2 points),
\item[(e)] one of the 2 points $\mathscr C_2\cap \mathscr C_5$
(consisting of 1 point).
\end{enumerate}
We show that if a hyperquadric $Q$ goes through these 12 points,
then $Q$ automatically contains all the five double curves.
In fact, from (a)--(c), $Q$ contains 5 points on $\mathscr C_1$,
and hence $\mathscr C_1\subset Q$.
If $Q$ further goes through the two points in (d),
$Q$ passes through 8 points on $\mathscr C_4$,
which means $\mathscr C_4\subset Q$.
Furthermore from the final point (e),
$Q$ goes through 5 points on $\mathscr C_2$, which means
$\mathscr C_2\subset Q$.
This implies that $Q$ passes through 8 points on $\mathscr C_5$,
and therefore $\mathscr C_5\subset Q$.
These mean that $Q$ contains 6 points on $\mathscr C_3$,
meaning $\mathscr C_3\subset Q$.
Thus the hyperquadric $Q$ contains all the five double curves.
Since $h^0(\mathscr O_{\CP^4}(2))=15$,
these hyperquadrics form a $2$-dimensional subsystem.
Any one of these $Q$-s which is different from $Y$ 
gives the required quadratic hypersurface.

The case of type III can be shown in a similar way.
We omit the detail.
\proofend

\subsection{Defining equation of the branch divisor}
We are ready for determining defining equation
of the  branch divisor of the anticanonical map.
As before let $\pi:\CP^4 \to \CP^2$ be a linear projection,
$\Lambda$ a conic in the target plane, and 
$Y\subset\CP^4$ the scroll over $\Lambda$. 
We know $Y$ is the image of the anticanonical map of the twistor spaces.
We fix any homogeneous coordinates 
$(z_0,z_1,z_2)$ on the above plane
such that the conic $\Lambda$ is defined by 
the equation $z_0^2 = z_1z_2$.

\begin{theorem}\label{thm:main}
Let $Z$ be  a Moishezon twistor space on $4\CP^2$
whose anticanonical map is (rationally) of degree two over the image,
and $\Phi:Z\to Y\subset\CP^4$ the anticanonical map.
Then the branch divisor of \,$\Phi$ is an intersection of the scroll \,$Y$ with a 
quartic hypersurface
defined by the following equation:
\begin{enumerate}
\item[
(i)] If $Z$ is of type II, the equation is of the form
\begin{align}\label{eqn:q1}
z_0z_1z_3z_4  = Q(z_0,z_1,z_2,z_3,z_4)^2,
\end{align}
where $Q$ is a quadratic polynomial such that the discriminant 
of the quadratic form $Q(0,0,z_2,z_3,z_4)$ is zero.

\item[(ii)]
If $Z$ is of type III, the equation of the hyperquartic is of the form
\begin{align}\label{eqn:q2}
z_0z_1z_4 (z_0-\lambda z_1) = 
Q(z_0,z_1,z_2,z_3,z_4)^2,
\end{align}
where $\lambda$ is a real number satisfying $\lambda\neq0,1$, and $Q$ is a quadratic polynomial such that the discriminant 
of the quadratic form $Q(0,0,z_2,z_3,z_4)$ is zero.

\end{enumerate}
\end{theorem}


%

\begin{remark}{\em
Before proceeding to the proof,
because of the coherency of the four types of the twistor spaces, we here write the equations 
of the branch divisors in the cases of 
type I and type IV. For the case of type I, under the above normalization
for the equation of the conic $\Lambda$, the equation of the hyperquartic is of the form
\begin{align}
z_0z_3 z_4  f (z_0,z_1,z_2,z_3,z_4) = Q (z_0,z_1,z_2,z_3,z_4) ^2,
\end{align}
where $f$ and $Q$ are  linear and quadratic forms respectively.
(No constraint for the discriminant of $Q$ for this case.)
For the case of type IV, under the same normalization, the equation is of the form
\begin{align}
z_0z_1(z_0-\lambda_1z_1)(z_0-\lambda_2 z_1) = Q (z_0,z_1,z_2,z_3,z_4) ^2,
\end{align}
where $\lambda_1$ and $\lambda_2$ are distinct real numbers satisfying 
$\lambda_1,\lambda_2\not\in\{0,1\}$.
Further  the discriminant 
of the quadratic form $Q(0,0,z_2,z_3,z_4)$ vanishes.

Although these might look complicated, a general principle is simple:
regardless of the types, the equation of the hyperquartic is always of the form
\begin{align}\label{gform}
{\text{(product of four linear polynomials)}} = Q (z_0,z_1,z_2,z_3,z_4) ^2,
\end{align}
and according to each of the degenerations
I $\to$ II $\to$ III $\to$  IV,
one of the linear polynomial degenerates from 
those with 5 variables $z_0,z_1,z_2,z_3,z_4$ to those with 3 variables 
$z_0,z_1,z_2$; 
geometrically the hyperplanes (in $\CP^4$) defined by the linear polynomials
degenerate from general ones to those which contain 
the singular line $l$ of the scroll $Y$.
So the `absolute value' of the type coincides with the number of 
the linear forms in the left-hand-side of \eqref{gform}
which belong to $\mathbb C[z_0,z_1,z_2]$.

Thus there is a strong similarity with the case of $3\CP^2$ obtained by Kreussler-Kurke \cite[p.\,49,50]{KK92}
(see also Poon \cite{P92}).
}
\end{remark} 

\vsp
\noindent
{\em Proof of Theorem \ref{thm:main}.}
For an algebraic subset $X\subset\CP^n$,
we denote by $I_X\subset\mathbb C[z_0,\cdots,z_n]$
the homogeneous ideal of $X$.
By Proposition \ref{prop:branch} there exists
a hyperquartic $\mathscr B$ in $\CP^4$ such that 
the branch divisor $B$ of the anticanonical map
is given as $Y\cap \mathscr B$.
Let $F=F(z_0,\cdots,z_4)$ be a defining quartic
polynomial of $\mathscr B$, 
where $z_0,\cdots,z_4$ are homogeneous coordinates
on $\CP^4$,
chosen in such a way that the projection $\pi:\CP^4\to \CP^2$
is given by $(z_0,\cdots,z_4)\mapsto (z_0,z_1,z_2)$.
In particular we have $l= \{z_0=z_1=z_2\}$.
Obviously $F$ is defined only up to an ideal
$I_{Y}\subset\mathbb C[z_0,\cdots,z_4]$.
Let $Q$ be a defining polynomial of the hyperquadrics
in $\CP^4$ containing all the double conics,
whose existence was proved in 
Proposition \ref{prop:Q}.
Let $\mathbb P_i\subset Y$
be the plane $\pi^{-1}(\lambda_i)$,
where $1\le i\le 3$ in the case of type II
and $1\le i\le 4$ in the case of type III.
Then as $(F|_{\mathbb P_i}) = 2\mathscr C_i
= (Q^2|_{\mathbb P_i})$
as divisors on the plane $\mathbb P_i$,
there exists a constant $c_i$ such that $F-c_iQ^2 \in
I_{\mathbb P_i}\subset\CC[z_0,\cdots,z_4]$.
If $c_i\neq c_j$ for some $i\neq j$,
we obtain $Q^2\in I_{\mathbb P_i} + I_{\mathbb P_j}$.
Further the last ideal is readily seen to be equal to
$ I_{\mathbb P_i\cap\mathbb P_j}$, and therefore
equals to $I_l=(z_0,z_1,z_2)\subset\CC[z_0,\cdots,z_4]$.
Hence $Q\in (z_0,z_1,z_2)$.
But this means that the divisor  $(Q|_{\mathbb P_i})$
contains $l$, which contradicts
the structure of the double conics (including the splitting one)
obtained in Section \ref{ss:dc}.
Therefore $c_i=c_j$ for any double conics
$\mathscr C_i$ and $\mathscr C_j$.

Next for the double quartic curve $\mathscr C_k$,
so that $k=4,5$ in the case of type II and
$k=5$ for the case of type III,
since $(F|_{H_k\cap Y}) = 2\mathscr C_k
= (Q^2|_{H_k\cap Y})$,
there exists a constant $c_k$ such that 
$F - c_kQ^2\in I_{H_k\cap Y}=(z_k) + I_Y$.
So taking a difference with $F-c_1Q^2\in I_{\mathbb P_1}$,
we obtain that $(c_1-c_k)Q^2 \in (z_k) + I_Y + 
I_{\mathbb P_1}$.
But since $\mathbb P_1\subset Y$, we have
$I_{\mathbb P_1}\supset I_Y$,
and therefore 
$(c_1-c_k)Q^2 \in (z_k) + I_{\mathbb P_1}$.
Hence if $c_1\neq c_k$ we have $Q^2\in (z_k) + I_{\mathbb P_1}$, which means $Q^2|_{\mathbb P_1}\in (z_k|_{\mathbb P_1})$. 
Since $k>2$, this means that the divisor $(Q^2)|_{\mathbb P_1}$ contains a line $(z_k)$ on 
the plane $\mathbb P_1$ as an irreducible component,
which  again contradicts the irreducibility of $\mathscr C_1$.
Therefore we have $c_1=c_k$ for any double 
quartic curve $\mathscr C_k$.
By rescaling we can suppose $c_i=1$ 
for any $1\le i\le 5$.
Thus  we have
\begin{align}
F - Q^2 &\in I_{\mathbb P_1}\cap 
I_{\mathbb P_2}\cap 
I_{\mathbb P_3}\cap ((z_3) + I_Y)
\cap  ((z_4) + I_Y)
&{\text{for the case of type II}},
\label{gaeof9}
\\
F - Q^2 &\in I_{\mathbb P_1}\cap 
I_{\mathbb P_2}\cap 
I_{\mathbb P_3}\cap 
I_{\mathbb P_4}
\cap  ((z_4) + I_Y)
&{\text{for the case of type III}}.
\end{align}

For the case of type III, 
by using the last ideal,
we can write
$ F- Q^2 = z_4 f+ g$, where 
$f$ is a cubic polynomial and 
$g$ is a quartic polynomial in $I_Y$.
This readily means $z_4f\in I_{\mathbb P_i}$
for any double conic $\mathscr C_i$
(namely for $i=1,2,3,4$).
For $i=1,2,3$,
let $l_i$ be the line going through the two points
$\lambda_i$ and $\lambda_4$, and $f_i$ 
a defining linear polynomial of $l_i$.
Then we have 
\begin{align}
Y \cap \pi^{-1} (l_i) = \mathbb P_i \cup \mathbb P_4,
\quad 1\le i\le 3,
\end{align}
and therefore
$$I_{\mathbb P_i}\cap I_{\mathbb P_4}
=I_{\mathbb P_i\cup\mathbb P_4} = I_{Y\cap  \pi^{-1}(l_i)}
= I_Y + I_{\pi^{-1}(l_i)}.$$
Hence from $z_4 f \in I_{\mathbb P_i}\cap I_{\mathbb P_4}$
 we can write $z_4 f = yg + f_ih$
where $y$ is a defining quadratic polynomial of the scroll $Y$.
If we write $g= z_4 g_1 + g_2$ and $h= z_4 h_1 + h_2 $
in a way that $g_2$ and $h_2$ do not involve $z_4$,
then we compute
$
y g + f_i h = z_4 ( yg_1 + f_i h_1) + 
(yg_2 + f_i h_2).
$
As this equals $z_4 f$ and
$yg_2 + f_i h_2$ does not involve $z_4$, we obtain $yg_2 + f_i h_2 =0$.
Hence, since $y$ and $f_i$ are mutually prime
from  irreducibility of the conic $\Lambda$, 
we can write $h_2 = y h_3$ by some quadratic polynomial $h_3$.
Thus we obtain
$$
z_4 f = z_4 f_i h_1 + y(g+ f_ih_3).
$$
Repeating a similar argument we can  pull out the 
linear polynomials $f_1, f_2$ and $f_3$ one by one, and 
 it follows that $z_4f$ can be written of 
the form $z_4 f_1f_2f_3+ \eta$, where $\eta\in I_Y$.
Therefore we obtain 
\begin{align}\label{final33}
F - Q^2 = z_4 f_1 f_2 f_3 + (g+\eta), \quad g+\eta\in I_Y.
\end{align}
By usual ${\rm PGL}(3,\CC)$-action 
we can normalize the homogeneous coordinates $(z_0,z_1,z_2)$
in such a way that 
$$
\lambda_4 = (0,0,1), \quad
\lambda_1 = (0,1,0), \quad
\lambda_2 = (1,0,0)
$$
hold.
Then we may suppose $f_1 = z_0$ and $f_2 = z_1$.
Moreover, as $l_3\ni \lambda_4$ and all $\lambda_i$-s are real and
mutually distinct, $f_3$ must be of the form 
of the form $z_0-\lambda z_1$, where $\lambda\in\RR\backslash\{0,1\}$.
Thus disposing $g+\eta$, \eqref{final33} means that
 a defining equation of the 
hyperquartic $\mathscr B$ can be taken in the form
\eqref{eqn:q2} in Theorem \ref{thm:main}.

For the case III, we still remain to show the claim about discriminant of $Q$.
As $\lambda_4  = (0,0,1)$ from the above choice of the coordinates,
the image $\Phi(S_4^+)$ and $\Phi(S_4^+)$ and over the plane
$\{z_0=z_1=0\}$.
Further as in Proposition \ref{prop:dsg84ga} (ii), 
these images are mutually distinct lines.
Hence $Q(0,0,z_2,z_3,z_4)$ must split into two linear forms,
as claimed.

For the case of type II,
in the above argument for the case of type III,
we replace the role of the fourth double conic
$\mathscr C_4$ by the third one $\mathscr C_3$,
and also the role of the remaining double conics
$\mathscr C_i$ $(i=1,2,3)$ by $\mathscr C_1$ and
$\mathscr C_2$.
This gives a similar expression 
\begin{align}\label{98gup}
F- Q^2 = z_3 f_1 f_2 f+ g,
\end{align}
where $f_i$ $(i=1,2)$ is a defining equation of 
the line connecting the points $\lambda_3$ and $\lambda_i$,  $f$ is a linear polynomial in $z_0,\cdots,z_4$, and $g\in I_Y$.
Further, from \eqref{gaeof9} we have 
$F-Q^2 \in (z_4) + I_Y$.
Hence we can write
$$
z_3 f_1 f_2 f = z_4 h + g_1, \quad
g_1\in I_Y.
$$
If we write $f= cz_4 + \zeta$ where $c\in\CC$ and $\zeta$ is a 
linear polynomial without the variable $z_4$,
we obtain $ z_4 (cz_3 f_1 f_2 - h) = -z_3 f_1 f_2 \zeta + g_1$.
Therefore, since $f_1f_2\zeta$ does not involve $z_4$, we have $\zeta=0$, and we obtain  $f=cz_4$.
By \eqref{98gup} this means 
$$
F- Q^2 = cz_3z_4f_1f_2 + g,
\quad g\in I_Y.
$$
If $c=0$, then $F-Q^2 \in I_Y$, which means 
$F|_Y= Q^2 |_Y$.
This means that the branch divisor of the 
double covering $\Phi_3:Z_3\to Y$ is non-reduced.
But of course this cannot happen since the restriction of $\Phi_3$ to a general fiber of $\pi:Y\to \Lambda$ is 
a non-reduced quartic curve
(Proposition \ref{prop:bian1}).
Therefore we have $c\neq 0$, and we obtain 
$F-Q^2 = z_3z_4 f_1 f_2 + g$ with $g\in I_Y$.
Then by the same argument in the case of type III,
we can normalize the homogeneous coordinates $(z_0,z_1,z_2)$
in a way that $f_1 = z_0$ and $f_2 = z_1$.
This means that a defining equation of the 
hyperquartic $\mathscr B$ can be taken in the form
\eqref{eqn:q1} in Theorem \ref{thm:main}.
The remaining claim about $Q$ follows in exactly the same way
as in the case of type III,
if we use Lemma \ref{lemma:anti2} (iii) instead of Proposition \ref{prop:dsg84ga} (ii).
\proofend

\section{Dimension of the moduli spaces}
\label{s:moduli}
In this section we compute dimension of the moduli space of the present twistor spaces.
For the case of type I, this was done in \cite[Section 5.1]{HonDS4_1},
but the argument in the paper does not work
in the cases of type II and type III.

Let $Z$ be a twistor space on $4\CP^2$ whose anticanonical map
is of degree two, and suppose that the type of $Z$ is 
I, II or III.
(We include type I since the present argument reproves 
a result in \cite{HonDS4_1}.)
Let $S\in |F|$ be a real irreducible member.
Then by a similar argument to \cite[Proposition 5.1]{HonDS4_1},
for the tangent sheaf of $Z$ we have
\begin{align}\label{hfgeg7}
H^i(\Theta_Z) = 0 {\text{ for }} i\neq 1,
\quad h^1(\Theta_Z)=13.
\end{align}
(Note that if $Z$ is of type IV, this is not the case.)
Also, it is not difficult to show 
\begin{gather}\label{73hew}
H^i(\Theta_S) = 0 {\text{ for }} i\neq 1,
\quad h^1(\Theta_S)=10,\\
h^0(K_S^{-1}) = 1, 
\quad H^i(K_S^{-1})=0{\text{ for }} i\neq 0.
\label{839l}\end{gather}
Let $\Theta_{Z,S}$ denote the subsheaf 
of  $\Theta_Z$ consisting 
of germs of vector fields which are tangent to $S$,
and write $\Theta_Z(-S) := \Theta_Z\otimes\mathscr O_Z(-S)$.
Then by various standard exact sequences of sheaves including these, and noting $N_{S/Z}\simeq K_S^{-1}$,
we obtain the following commutative diagram of cohomology groups
on $Z$ and $S$:
$$
\begin{CD}
@. 0 @. 0 @. 0 @.\\
@. @VVV @VVV @VVV @.\\
0 @>>> H^0(\Theta_Z|_S)  @>>> H^0(K_S^{-1})
@>>> H^1(\Theta_S) @>>> H^1(\Theta_Z|_S)\\
@. @VVV @VV{\delta}V @| @.\\
0 @>>> H^1(\Theta_Z(-S)) @>>> H^1(\Theta_{Z,S})
@>{\alpha}>> H^1(\Theta_S) @>>> 0\\
@. @VVV @VV{\beta}V @VVV @.\\
0 @>>> H^1(\Theta_Z) @= H^1(\Theta_Z) @>>> 0\\
@. @. @VVV\\
@. @. 0
\end{CD}
$$
From the middle column of this diagram we obtain 
$h^1(\Theta_{Z,S}) = 14$ by
\eqref{hfgeg7} and \eqref{839l},
which means $h^1(\Theta_Z(-S)) = 4$ from 
the middle row and \eqref{73hew}.
In particular, the Kuranishi family of
deformations of the pair $(Z,S)$ is 14-dimensional.

In order to compute the dimension of the moduli spaces,
we first recall that our twistor spaces
(of types I--III) can be 
characterized by the property that they have 
a rational surface $S$ with particular structure 
as a member of the system $|F|$.
For type II twistor spaces the structure of $S$ is 
described in the proof of Proposition \ref{prop:reducible} (i) in terms of blownup points of the birational morphism
$\epsilon:S\to \CP^1\times\CP^1$.
In particular, all freedom for $S$
that can be contained in the twistor spaces 
of type II is to move four points on a fixed cycle of anticanonical curve on $\CP^1\times\CP^1$
(the remaining four points does not contribute
for deforming complex structure),
and therefore they constitute 4-dimensional family
of deformations of $S$.
Via the Kodaira-Spencer map, this family determines
a 4-dimensional subspace of $H^1(\Theta_S)$,
for which we denote by $V$.

We have $\dim \alpha^{-1}(V) = \dim V + h^1(\Theta_Z(-S))= 8$.
The tangent space of the moduli space of 
twistor spaces of type II can be considered as
the space $\beta(\alpha^{-1}(V))
\subset H^1(\Theta_Z)$.
The image of the map $\delta$ in the diagram
corresponds to deformations of $(Z,S)$ 
that can be obtained by moving $S$ in $Z$,
and of course they do not give a non-trivial deformation
of $Z$.
But from the characterization of $Z$ by
the complex structure of $S$, even if we move $S$ in $Z$, its complex structure
cannot go away from the above 4-dimensional family of $S$.
This means that the image
of $\delta$ is contained in $\aaa^{-1}(V)$.
Thus the tangent space  of the moduli space of 
twistor spaces of type II can be identified with
the quotient space
\begin{align}\label{8fds82}
\aaa^{-1}(V)/\delta H^0(K_S^{-1}),
\end{align}
and this is 7-dimensional by \eqref{839l}.

For the case of type III, the same argument works 
except that the subspace
$V\subset H^1(\Theta_S)$ 
becomes 2-dimensional in this case.
Consequently the tangent space of the moduli space
is again identified with the quotient space
\eqref{8fds82}, which is 5-dimensional this time.

This way we conclude that the moduli space of 
twistor spaces  is 7-dimensional for the case of 
type II,
and $5$-dimensional  for the case of  type III.
We note that as obtained in \cite[Section 5.1]{HonDS4_1},
for the case of type I the moduli space is $9$-dimensional.
This can also be seen from the above argument if we notice
$\dim V=6$ for the case of type I.
Thus according to the degenerations
I $\to$ II and II $\to$ III,
the dimension of the moduli spaces drops by two.
On the other hand for the case of type IV
the moduli space is $4$-dimensional
as obtained in \cite{HonDSn1}.
This discrepancy comes from the fact that 
twistor spaces of type IV admit a
non-trivial $\CC^*$-action,
while other three types of spaces do not.

\begin{remark}
{\em
Looking the equations of the quartic hypersurfaces in 
Theorem \ref{thm:main}, one may think that the dimension of the moduli spaces
computed above contradicts Theorem \ref{thm:main},
because the number of parameters involved in the equation \eqref{eqn:q2}
in the case of type III is greater than 
those for the equation \eqref{eqn:q1} in the case of type II. 
But it is not correct, since 
the elements in ${\rm{PGL}}(5,\CC)$ preserving 
the form of \eqref{eqn:q2} is greater than those for \eqref{eqn:q1}.
}
\end{remark}

%

\end{document}